\documentclass[12pt]{amsart}

\newcommand\real{{\mathrm I}\!{\mathrm R} }

\newcommand\rat{{\mathrm Q}\kern-.65em {}^{{}_/ }}
\newcommand\intg{{{\mathrm Z}\mkern-5.5mu{\mathrm Z}}}

\newtheorem{remarkk}{Remark (question)}

\newtheorem{theorem}{Theorem}
\newtheorem{lemma}{Lemma}
\newtheorem{proposition}{Proposition}

\begin{document}
\title{On the Growth of hyperbolic geodesics in rank 1 manifolds}
\author{Abdelhamid Amroun}
\date{}
\maketitle
\address{\begin{center}
Universit\'e Paris-Sud, D\'epartement de Math\'ematiques,
CNRS UMR 8628, 91405 Orsay Cedex France
\end{center}}

\begin{abstract} We give a formula for the topological pressure of the
geodesic flow of a compact rank 1 manifold in terms of the growth of
the number of closed hyperbolic (rank 1) geodesics. We derive 
an equidistribution result for these
geodesics with respect to equilibrium states. This generalize
partially a result of G. Knieper \cite{kni} to non constant potentiels.
\end{abstract}
\section{Introduction and main result}
By a result of G Knieper \cite{kni}, we know that the geodesic flow on
a compact rank 1 manifold $M$ admits a unique measure of 
maximal entropy, concentrated on the open and dense set of
regular vectors in $T^{1}M$.
Knieper proved that this measure is
approximated by probability measures supported on a finite number of
regular closed geodesics.
This result generalize previous well known results for negatively curved
manifolds  \cite{ano} \cite{bow} \cite{bowe} \cite{led} \cite{mar}
\cite{marg} \cite{cby}.

In this paper we  consider the case of non constant potentials.
We obtain a formula expressing the topological
pressure as an exponential growth of the number of weighted regular closed
geodesics representing different free homotopy classes.
As a consequence we give an equidistribution result for weighted
closed regular geodesics to an equilibrium state. These results extend and
strengthens previous one by G Knieper (\cite{kni} Proposition $6.4$)
and  M Pollicott \cite{pol}. 
The proof uses the Anosov's closing lemma for compact
manifolds of nonpositive curvature \cite{BBS} and also a
Riemannian formula for the topological pressure by G P Paternain \cite{pat}.  

Let $M=X/\Gamma$ be a compact Riemannian manifold
of nonpositive curvature where $X$ is the universal cover and $\Gamma$
is the group of deck transformations of $X$.
The rank of a vector $v\in T^{1}M$ is the dimension of
the space of all parallel Jacobi fields along the geodesic defined by
$v$. The rank of the manifold is the minimal rank of all tangeant
vectors. We will assume that $M$ is a rank
1 manifold (this includes manifolds of negative sectional curvature
where all the geodesics are of rank 1). In fact,
by a rigidity result of Ballmann
\cite{bal} and Burns-Spatzier \cite{bs} ``most of'' compact manifolds of
nonpositive curvature are rank 1.
By a regular vector (resp regular geodesic) we will mean a rank
1 vector (resp a geodesic defined by a rank 1 vector). A geodesic is
called hyperbolic if it is regular, extending thus the notion of
hyperbolicity to rank $1$ manifolds.  
Let $\mathcal{R}_{reg}$ be the open subset of $T^{1}M$ of regular vectors. It
is dense in $T^{1}M$ if $M$ is of finite volume \cite{ball}.
Let $\nu$ be the Knieper's measure of maximal entropy of the
geodesic flow of the rank 1 maniflod $M$. We have
$\nu(\mathcal{R}_{reg})=1$ and the
complement $\mathcal{S}_{ing}$ of $\mathcal{R}_{reg}$ is an invariant
closed subset of the unit tangent bundle. The growth of closed
geodesics in the ``singular part'' $\mathcal{S}_{ing}$ can be
exponential \cite{gromov} as well as subexponential \cite{kni}. 
In this paper we concentrate on the regular set, but it will be
interesting to investigate the $\mathcal{S}_{ing}$-part. 

Two elements $\alpha, \beta \in \Gamma$ are equivalent if
and only if there exists $n, m \in \intg$ and $\gamma \in \Gamma$ such
that $\alpha^{n}=\gamma \beta^{m} \gamma^{-1}$. 
Denote by $[\Gamma]$ the set of equivalence classes of elements in
$\Gamma$.
Classes in $[\Gamma]$
are represented by elements in $\Gamma$ which have a least period
(primitive elements):
\[
[\alpha]=\{\gamma \alpha_{0}^{m} \gamma^{-1}: \alpha_{0}\in \Gamma,
\ \alpha_{0} \ primitive, \  \gamma \in \Gamma \}.
\]
Let $x_{\alpha}$ be the point in $X$ such that $d(x_{\alpha}, \alpha
x_{\alpha})=\inf_{p\in X}d(p,\alpha p)$ ($M$ is compact). The axis
trougth $x_{\alpha}$ and $\alpha x_{\alpha}$ projects onto
a closed geodesic in $M$ with prime period
$d(x_{\alpha}, \alpha x_{\alpha}):=l(\alpha_{0})$.
We set $l([\alpha]):=l(\alpha_{0})$, i.e
\[
l([\alpha])=\min \{l(\gamma): \gamma \in [\alpha]\}=l(\alpha_{0}).
\]
We will denote by
$\Gamma_{hyp}\subset \Gamma$ the subset of those elements with hyperbolic
axis. Then $[\Gamma_{hyp}]$ is the set of conjugacy classes representing 
geometrically distinct hyperbolic closed geodesics. Finally given a
function $f$ on $T^{1}M$ the notation $\int_{[\alpha]}f$, $[\alpha] \in
[\Gamma_{hyp}]$, stands for the integral of $f$ along the unique closed
geodesic representing the class $[\alpha]$. If this geodesic is given by 
 $\phi_{s}v_{[\alpha]}, 0\leq s \leq l([\alpha])$ for some
$v_{[\alpha]}\in T^{1}M$, then,  
\[
\int_{[\alpha]}f:=\int_{0}^{l([\alpha])}f(\phi_{s}v_{[\alpha]})ds:=
\delta_{[\alpha]} (f). 
\]
Given a continuous function $f$ on $T^{1}M$, let $\mu_{t}:=\mu_{t}^{f}$ 
be the flow invariant probability measures  supported on a finite number of
 hyperbolic closed geodesics defined on continuous functions $\omega$ by,
\[
\mu_{t}(\omega):=
\frac{ \sum_{([\alpha] \in \Gamma_{hyp}:l([\alpha])\leq
  t)}e^{\int_{\alpha}f}\delta_{[\alpha]}(\omega)} 
{ \sum_{([\alpha] \in \in \Gamma_{hyp}: l([\alpha])\leq
  t)}e^{\int_{\alpha}f}}.
\]
Here is the main result of the paper.
\begin{theorem} Let $M=X/\Gamma$ be a compact rank 1 manifold equipped with a
  $C^{\infty}$ Riemannian metric and $f\in
  \mathcal{C}_{\real}(T^{1}M)$. Then
\begin{enumerate}
\item 
\begin{equation*}
\lim_{t\rightarrow +\infty}\frac{1}{t}
\log \sum_{[\alpha] \in  [\Gamma_{hyp}]: l([\alpha])\leq
t}e^{\int_{\alpha}f}=P(f).
\end{equation*}
\item The accumulation points of $\{\mu_{t}\}$ with respect to the 
topology of weak convergence of measures, are equilibrium states of
the geodesic flow corresponding to 
the potential $f$. Moreover, for any
open neighborhood $V$ in $\mathcal{P}(T^{1}M)$ of the subset of
equilibrium states $\mathcal{P}_{e}(\phi)$ we have,
\[
\lim_{t\rightarrow +\infty}
\frac{ \sum_{([\alpha] \in \Gamma_{hyp}:l([\alpha])\leq
  t, \ \delta_{[\alpha]} \in V)}e^{\int_{\alpha}f}} 
{ \sum_{([\alpha] \in \in \Gamma_{hyp}: l([\alpha])\leq
  t)}e^{\int_{[\alpha]}f}}=1,
\]where the convergence is exponential.
\end{enumerate}
\end{theorem}
In the part $(2)$ of Theorem $1$, the condition $C^{\infty}$ on the
metric is necessary since we need the upper-semicontinuity of the
entropy map \cite{new}. 

Let $d$ be the distance on $T^{1}M$ induced by the
Riemannian metric of $M$.
Consider the metric $d_{t}$ on $T^{1}M$, defined for all $t>0$ by
\[
d_{t}(u,v):=\sup_{0\leq s \leq t}d(\phi^{s}(u), \phi^{s}(v)).
\]
Following \cite{BBS} we denote by $P(t,\epsilon)$ the maximal number
of regular vectors $v\in T^{1}M$ which are $\epsilon$-separated in the metric
$d_{t}$ and for which $\phi^{t(v)}v=v$ for some $t(v) \in [t,
t+\epsilon]$. Let $E(t,\epsilon)$ be the set defined above with $\#
E(t,\epsilon)= P(t,\epsilon)$.
 
The following is the Lemma 5.6 from \cite{BBS} for rank 1 manifolds
and continuous potentials (see Lemma 1 $(5)$ below). 
\begin{proposition} 
Set
$\int_{c_{v}}f:=\int_{0}^{l(c_{v})}f(\phi^{t}(v))dt$, where $c_{v}$ is
the closed geodesic defined by $v\in  E(t,\epsilon)$ and $l(c_{v})$ is the
period of $v$. Then, 
\[
\lim_{\epsilon \rightarrow 0}\liminf_{t\rightarrow +\infty}\frac{1}{t}
\log \sum_{v\in E(t,\epsilon)}e^{\int_{c_{v}}f}=
\lim_{\epsilon \rightarrow 0}\limsup_{t\rightarrow +\infty}\frac{1}{t}
\log \sum_{v\in E(t,\epsilon)}e^{\int_{c_{v}}f}=P(f).
\]
\end{proposition} 
\section{Proofs}
\subsection{Topological pressure}
We recall the notion of topological pressure \cite{wal}.
Let $t>0$ and $\epsilon >0$. A subset $E\subset T^{1}M$ is a
$(t,\epsilon)$-separated if $d_{t}(u,v)>\epsilon$ for $u \ne v \in E$.
Set
\[
r(f; t,\epsilon):=\sup_{E}\sum_{\theta \in
E}e^{\int_{0}^{t}f(\phi_{t}(\theta))dt} 
\] 
where $\sup$ is over all $(t,\epsilon)$-separated subsets $E$; and 
\[
r(f;\epsilon):=\limsup_{t\rightarrow \infty}\frac{1}{t}
\log r(f; t,\epsilon).
\]
Then the topological pressure of the geodesic flow corresponding to the
potential de $f$ is the number, 
\begin{equation}
P(f)=\lim_{\epsilon \rightarrow 0}r(f;\epsilon).
\end{equation}
The topological entropy $h_{top}$ is $h_{top}=P(0)$.
We denote by $\mathcal{P}(T^{1}M)$ the set of all probability mesures on
$T^{1}M$ with the weak topology of mesures, and let $\mathcal{P}(\phi)$ be
the subset of invariant probability mesures of the flow.
The entropy of a probability measure $m$ is denoted $h(\mu)$ \cite{wal}.
All these objects satisfy the following variational principle \cite{wal}
\begin{equation}
P(f)=\sup_{\mu\in \mathcal{P}(\phi)}\left (h(\mu)+\int_{T^{1}M}fd\mu
\right ).
\end{equation}
An equilibrium state $\mu_{f}$ satisfies,
\begin{equation}
h(\mu_{f})+\int_{T^{1}M}fd\mu_{f}=P(f).
\end{equation}
When the Riemannian metric of the manifold $M$ is $C^{\infty}$ then by
a result of Newhouse \cite{new} the entropy map $m\rightarrow h(m)$ is
upper semicontinuous and then $h_{top}<\infty$. Consequently, the
set $\mathcal{P}_{e}(f)$ of equilibrium states
is a non empty closed and convex subset of $\mathcal{P}(\phi)$.
\subsection{Proof of Theorem $1$ $(1)$}
Let $\nu$ be the Knieper's measure. Let $N(t, \epsilon, 1-\delta,
\nu)$ be the minimal number of $\epsilon$-balls in the metric $d_{t}$
which cover a set of measure at least $1-\delta$. Since $\nu$ is
the unique measure of maximal entropy, we can apply Lemma 5.6 in
\cite{BBS} to this measure. 
\begin{lemma}[\cite{BBS}] There exists $\delta >0$ such that for all
$\epsilon >0$, there exists $t_{1}>0$ such that 
\begin{equation}
P(t,\epsilon)\geq N(t, \epsilon, 1-\delta, \nu)
\end{equation}
for any $t\geq t_{1}$. In particular, we have
\begin{equation}
\lim_{\epsilon \rightarrow 0}\liminf_{t\rightarrow +\infty}\frac{\ln
  P(t,\epsilon)}{t}= \lim_{\epsilon \rightarrow 0}\limsup_{t\rightarrow
  +\infty}\frac{\ln P(t,\epsilon)}{t}=h_{top}.
\end{equation}
\end{lemma}
We fix $\epsilon >0$. Recall that $E(t,\epsilon)$ is the maximal set
defined above with $\# E(t,\epsilon)=P(t,\epsilon)$.
From Lemma 1 $(4)$, if $N(t,\epsilon, 1-\delta, \nu)$ is the minimal
number of $\epsilon$-balls $B_{t}(v_{i},\epsilon)$ in the metric $d_{t}$,
which cover the whole space $T^{1}M$, then
$P(t,\epsilon)=N(t,\epsilon, 1-\delta,\nu)$ for $t$ sufficiently large
(since each $\epsilon$-ball $B_{t}(v_{i},\epsilon)$ contains a unique
point from $E(t,\epsilon)$).  

Now, it suffices to prove Theorem $1$ ($1$) for Lipschitz functions $f$.
Suppose then $f$ Lipschitz and let $lip(f)$ be it's Lipschitz constant.
Consider a $(2\epsilon,t)$-separated set $E_{1}$ in
$T^{1}M$. Thus, two distinct vectors in $E_{1}$ lies in two distinct
$\epsilon$-balls above, so that $\# E_{1}\leq P(t,\epsilon)$. 
For each $\theta \in E_{1}$, we associate the unique point $v_{\theta} \in
E(t, \epsilon)$ such that $d_{t}(\theta,v_{\theta} )\leq
\epsilon$. Let $\tau_{\theta}$ the regular closed geodesic
corresponding to the regular periodic vector $v_{\theta}$
($\dot{\tau}_{\theta}(0)=v_{\theta}$), with period $l(\tau_{\theta})\in [t,
  t+\epsilon]$, and  
$[\tau_{\theta}]$ the corresponding free homotopy class.
There exists a constant $C>lip(f)$ such that,
\begin{eqnarray*}
&&\sum_{\theta \in E_{1}}e^{\int_{0}^{t} f(\phi_{s}(\theta))ds}\\
&\leq & e^{lip(f)\epsilon t} \sum_{[\tau_{\theta}]: \theta
\in E_{1},  t<l(\tau_{\theta})\leq t+\epsilon}
e^{\int_{0}^{t} f(\tau_{\theta}(s))ds}\\
&\leq & e^{C\epsilon t} \sum_{[\tau]: t<l(\tau)\leq t+\epsilon} 
e^{\int_{0}^{l(\tau)} f(\tau(s))ds} \\
&\leq & e^{C \epsilon t} \sum_{[\tau]: l(\tau)\leq t+1} 
e^{\int_{0}^{l(\tau)} f(\tau(s))ds}, 
\end{eqnarray*} 
where the sum is over all the hyperbolic closed geodesics which represent
different free homotopy classes and prescribed lengh.
Thus for all $\epsilon <\epsilon_{0}$ we obtain that
\begin{eqnarray*}
P(f;2\epsilon)&=& \limsup_{t\rightarrow \infty}\frac{1}{t}
\log \sup_{E_{1}}\left (\sum_{\theta \in
  E_{1}}e^{\int_{0}^{t} f(\phi_{s}(\theta))ds} \right )\\
&\leq& C\epsilon + \limsup_{t\rightarrow \infty}\frac{1}{t}
\log \left (\sum_{[\alpha] \in [\Gamma_{hyp}]:l([\alpha])\leq t}
e^{\int_{[\alpha]} f}\right ).
\end{eqnarray*}
We let $\epsilon \rightarrow 0$ gives,
\begin{eqnarray*}
P(f)&:=&\lim_{\epsilon \rightarrow 0}P(f;2\epsilon)\\
&\leq& \limsup_{t\rightarrow \infty}\frac{1}{t} 
\log \left (\sum_{[\alpha] \in [\Gamma_{hyp}]:l([\alpha])\leq t}
e^{\int_{[\alpha]} f}\right ).
\end{eqnarray*}
To show the reverse inequality, $\limsup \leq P(f)$, it suffices to
observe that the set  
of hyperbolic closed geodesics with length $\leq t$ and which represent
different free homotopy classes, is $\epsilon$-separated for all
$\epsilon < inj(M)$. 

We will show now that the $\liminf$ is bounded below by the
pressure. For this we use the following result \cite{pat}.
\begin{theorem}[G P Paternain \cite{pat}] Let $M$ be a closed connected
Riemannian manifold. 
If the metric of $M$ is of class $C^{3}$ and does not have conjugate
points, then for any $\delta >0$ we have,
\begin{equation*}
P(f)=\lim_{t\rightarrow \infty}
\frac{1}{t}
\log \int_{M\times M}
\left (\sum_{\gamma_{xy}: t-\delta <l(\gamma_{xy})\leq t}
e^{\int_{0}^{l(\gamma_{xy})}f(\gamma_{xy}(t),\dot{\gamma}_{xy}(t))dt}
\right )dxdy.  
\end{equation*}
\end{theorem}
For each $\delta >0$ and $(x,y)\in M\times M$ we consider the subset
of $T^{1}M$ defined by,
\[
E_{xy}:=\{\dot{\gamma}_{xy}(0): t-\delta < l(\gamma_{xy}) \leq t\}.
\]
It is finite for almost all $(x,y)\in M\times M$ \cite{patt}. As
consequence of the nonpositive curvature, the rank 1 manifold $M$ has no
conjugate points.
Thus, there exists a constant $\epsilon_{0}$, depending only on $M$,
such that $E_{xy}$ is $(2\epsilon,t)$-separated for $\epsilon
<\epsilon_{0}$. To see this, it suffices to lift every thing to the
universal cover of $M$ and use (\cite{pat} p135) and (\cite{has} p375). 
Thus, the preceeding arguments applied to $ E_{1}=E_{xy}$ give,
\begin{eqnarray*}
&& \liminf_{t\rightarrow \infty}\frac{1}{t} 
\log \int_{M\times M}\left ( \sum_{\theta \in E_{xy}}e^{\int_{0}^{t}
  f(\phi^{s}(\theta))ds}
\right )dxdy\\
&\leq& \liminf_{t\rightarrow \infty}\frac{1}{t} 
\log \left (\sum_{[\alpha] \in \Gamma_{hyp}:l([\alpha])\leq t}
e^{\int_{[\alpha]} f}\right ).
\end{eqnarray*}
But by Theorem 2, the left hand side of this inequality is a limit and
is equal to $P(f)$, which completes the proof.
\begin{remarkk}
Find a proof which did not appeal to Paternain's formula in Theorem 2 !
\end{remarkk}
\subsection{Proof of Proposition 2}
The proof of Proposition 2 follows from the above arguments. 
\subsection{Proof of Theorem $1$ $(2)$}
Consider the following functional which measures the
``distance'' of an invariant measure $m$ to the set of equilibrium states,
\[
\rho(m)=P(f)-\left (h(m)+\int fdm \right ).
\]
Set $\rho(E):=\inf (\rho(m): m\in E)$ for $E \subset
\mathcal{P}(\phi)$ and,
\[
[\Gamma_{hyp}](t):=\{[\alpha] \in [\Gamma_{hyp}]: l([\alpha])\leq t\}.
\]
\begin{lemma}  Let $M$ be a compact smooth rank 1 manifold and
  $f$ a continuous potential on $T^{1}M$.
Then, for any closed subset $K$ of $\mathcal{P}(\phi)$ we have,
\[
\limsup_{t\rightarrow +\infty }\frac{1}{t}\log
 \frac{\sum_{([\alpha] \in [\Gamma_{hyp}](t): \delta_{[\alpha]} \in K)}
 e^{\int_{[\alpha]} f}} 
{\sum_{[\alpha] \in [\Gamma_{hyp}](t)}
 e^{\int_{[\alpha]} f}} 
\leq  -\rho(K).
\]
\end{lemma}
We leave the proof of this lemma for later and finish the proof of
Theorem $1$. First, let $V$ be an open neighborhood of
$\mathcal{P}_{e}(f)$ and set $K=\mathcal{P}(\phi)\backslash V$. The
set $K$ is compact and $\rho(K)>0$. 
For $t$ sufficiently large we have by Lemma $3$, 
\[
1\geq  \frac{\sum_{([\alpha] \in [\Gamma_{hyp}](t): \delta_{[\alpha]} \in V)}
 e^{\int_{[\alpha]} f}} 
{\sum_{[\alpha] \in [\Gamma_{hyp}](t)}
 e^{\int_{[\alpha]} f}}\geq 1-e^{-t\rho(K)}.
\]
This proves the second assertion in part $(2)$ of the Theorem $1$. 

\subsubsection{Accumulation measures of $\mu_{t}$}
The proof of the fact that the accumulation measures of $\mu_{t}$ are in
$\mathcal{P}_{e}(\phi)$ follows \cite{am}.

We endow $\mathcal{P}(T^{1}M)$ with a distance $d$ compatible with the weak
star topology: take a countable base $\{g_{1}, g_{2}, \cdots\}$ of the
separable space 
$C_{\real}(T^{1}M)$, where $\|g_{k}\|=1$ for all $k$, and set:
\[
d(m, m'):=\sum_{k=1}^{\infty}2^{-k}\left |\int g_{k}dm - \int g_{k}dm'
\right |.
\] 
 
Let $V\subset \mathcal{P}(\phi)$ be a convex open 
neighborhood of $\mathcal{P}_{e}(f)$ and $\epsilon >0$.
We consider a finite open cover $(B_{i}(\epsilon))_{i\leq N}$ of
$\mathcal{P}_{e}(f)$ by balls of diameter $\epsilon$ all contained in $V$. 
Decompose the set $U:=\cup_{i=1}^{N}B_{i}(\epsilon)$ as follows,
\[
U=\cup_{j=1}^{N'}U_{j}^{\epsilon},
\]where the sets $U_{j}^{\epsilon}$ are disjoints (not necessarily
open ) and contained in one
of the balls $(B_{i}(\epsilon))_{i\leq N}$. We have
\[
\mathcal{P}_{e}(f) \subset U \subset V.
\] 
We fix in each $U_{j}^{\epsilon}$ an invariant probability measure $m_{j}$,
$j\leq N'$, and let $m_{0}$ be an invariant probability  measure
distinct from the above ones; for example take $m_{0} \in V\backslash U$.
Set for convenience,
\begin{equation}
\nu_{t}(E):= \frac{\sum_{([\alpha] \in [\Gamma_{hyp}](t):
 \delta_{[\alpha]} \in E)} e^{\int_{[\alpha]} f}} 
{\sum_{[\alpha] \in [\Gamma_{hyp}](t)}
 e^{\int_{[\alpha]} f}},\
 E\subset \mathcal{P}(\phi)
\end{equation}
and define,
\begin{equation}
\beta_{t}=
\sum_{j=1}^{N'}\nu_{t}(U_{j}^{\epsilon})m_{j}+
(1-\nu_{t}(U))m_{0}. 
\end{equation}
We have $\sum_{j=1}^{N'}\nu_{t}(U_{j}^{\epsilon})=\nu_{t}(U)$.
The probability measure $\beta_{t}$ lies in $V$ since it is a convex
combination of elements in the convex set $V$. Thus
\[
d(\mu_{t}, V)\leq d(\mu_{t},\beta_{t}). 
\]
We are going to show that
\[
d(\mu_{t},\beta_{t})\leq \epsilon \nu_{t}(U)+3\nu_{t}(U^{c}),
\]
where $U^{c}=\mathcal{P}(SM)\backslash U$.

Consider the measures $\mu_{t,V}$ on $SM$ defined by,
\[
\mu_{t,V}:=\frac{\sum_{([\alpha] \in [\Gamma_{hyp}](t):
 \delta_{[\alpha]} \in V)} 
 e^{\int_{[\alpha]} f}\delta_{[\alpha]}} 
{\sum_{[\alpha]\in [\Gamma_{hyp}](t)}
 e^{\int_{[\alpha]} f}}.
\]
By definition of $\mu_{t}$ and $\mu_{t,V}$ and
the fact that $U\subset V$,
\[
\sum_{k\geq 1}2^{-k}|\mu_{t}(g_{k})-\mu_{t,V}(g_{k})|\leq
\nu_{t}(U^{c}).
\]
It remains to show that 
\[
\sum_{k\geq 1}2^{-k}|\mu_{t,V}(g_{k})-\beta_{t}(g_{k})| \leq \epsilon
\nu_{t}(U)+\nu_{t}(U^{c}).
\] 
We have for all $k\geq 1$,
\[
|\mu_{t,V}(g_{k})-\beta_{t}(g_{k})|\leq A+B+C
\]where,
\[
A=
\frac{\sum_{j=1}^{N'}\sum_{[\alpha] \in [\Gamma_{hyp}](t):
\delta_{[\alpha]}\in U_{j}^{\epsilon}}
e^{\int_{[\alpha]}f}|\delta_{[\alpha]}(g_{k})-m_{j}(g_{k})|}
{\sum_{[\alpha] \in [\Gamma_{hyp}](t)}e^{\int_{[\alpha}] f}},
\]
\[
B=\frac{\sum_{[\alpha] \in [\Gamma_{hyp}](t):
\delta_{[\alpha]}\in V\backslash
U}e^{\int_{[\alpha]}f}\delta_{[\alpha]}(g_{k})} 
{\sum_{[\alpha] \in [\Gamma_{hyp}](t)}e^{\int_{[\alpha}] f}},
\]
\[
C=|(1-\nu_{t}(U))m_{0}(g_{k})|.
\]
Thus, since we have for all $k\geq 1$, $\|g_{k}\|=1$,  by definition
of $\nu_{t}$ we get,
\begin{eqnarray*}
&&\sum_{k\geq 1}2^{-k}|\mu_{t,V}(g_{k})-\beta_{t}(g_{k})|\\
&\leq& \epsilon
  \sum_{j=1}^{N'}\nu_{t}(U_{j}^{\epsilon})+\nu_{t}(U^{c})+
(1-\nu_{t}(U))\\
&=& \epsilon \nu_{t}(U)+\nu_{t}(U^{c}).
\end{eqnarray*}
Finally we have obtained that
\[
d(\mu_{t},\beta_{t})\leq \epsilon
\nu_{t}(U)+\nu_{t}(U^{c}). 
\]
This implies the desired inequality,
\[
d(\mu_{t},V)\leq \epsilon
\nu_{t}(U)+3\nu_{t}(U^{c}). 
\]
The set $U^{c}$ is closed, so we have $\lim_{t\rightarrow
\infty}\nu_{t}(U)=1$. Since $\epsilon$ is arbitrary, we
conclude that $\limsup_{t\rightarrow
\infty}d(\mu_{t},V)=0$. The neighborhood $V$ of
$\mathcal{P}_{e}(f)$ being arbitrary, this implies that all limit
measures of $\mu_{t}$ are contained in $\mathcal{P}_{e}(f)$. In
particular, if $\mathcal{P}_{e}(f)$ is reduced to one measure $\mu$, this
shows that $\mu_{t}$ converges to $\mu$.

\subsection{Proof of Lemma $3$}
We follow \cite{pol}.
The functional $\rho$ is lower semicontinuous (since $h$ is upper
semicontinuous) and $\rho \geq 0$. 
Set for any continuous function $\omega$ on $T^{1}M$,
\begin{equation}
Q_{f}(\omega):=P(f+\omega)-P(f).
\end{equation}
The fact that $Q_{f}$ is a convex and continuous is a consequence of
the same properties for $P$.
Using the variational principle, it is not difficult to see that
\[
Q_{f}(\omega)=\sup_{\mu \in \mathcal{P}(\phi)}\left (\int \omega
  d\mu-\rho(\mu)\right ).
\]
By duality we have for any invariant probabilité measure $m$,
\[
\rho(m)=\sup_{\omega \ continuous}\left (\int \omega
  dm-Q_{f}(\omega)\right ).
\]
With the notations introduced above, we have to prove that
\[
\limsup_{t\rightarrow +\infty}\frac{1}{t}\log \nu_{t}(K)\leq -\rho(K).
\]
Let $\epsilon >0$. There exists a finite number of continuous
functions $\omega_{1}, \cdots, \omega_{l}$ such that $K\subset
\cup_{i=1}^{l}K_{i}$, where 
\[
K_{i}=\{m\in \mathcal{P}(\phi): \int
\omega_{i}dm-Q(\omega_{i})>\rho(K)-\epsilon\}. 
\]
We have $\nu_{t}(K) \leq \sum_{i=1}^{l}\nu_{t}(K_{i})$ where
\[
\nu_{t}(K_{i})=\frac{\sum_{([\alpha] \in [\Gamma_{hyp}](t):
    \delta_{[\alpha}] \in K_{i})}  e^{\int_{[\alpha]} f}}
{\sum_{[\alpha] \in [\Gamma_{hyp}](t)}e^{\int_{[\alpha]} f}}.
\]
Then,
\begin{eqnarray*}
&&\sum_{[\alpha] \in [\Gamma_{hyp}](t):
\delta_{[\alpha]} \in K_{i}}  e^{\int_{[\alpha]} f}\\
&\leq&\sum_{[\alpha] \in [\Gamma_{hyp}](t):
\delta_{[\alpha]} \in K_{i}}  e^{\int_{[\alpha]} f}
e^{l([\alpha])(\int \omega_{i}d\delta_{[\alpha]}-
Q(\omega_{i})-(\rho(K)-\epsilon))}.
\end{eqnarray*}
Set $C:=\sum_{i\leq l}\sup(1, e^{-\delta(-Q(\omega_{i})-(\rho(K)-\epsilon))})$.
Thus, by taking into account the sign of $-Q(\omega_{i})-(\rho(K)-\epsilon)$,
\begin{eqnarray*}
&& \sum_{\alpha \in [\Gamma_{hyp}](t):
\delta_{[\alpha]} \in K_{i}}  e^{\int_{[\alpha]} f}\\
&\leq& \sum_{[\alpha] \in [\Gamma_{hyp}](t):
\delta_{[\alpha]} \in K_{i}}  e^{\int_{[\alpha]} (f+\omega_{i})}
e^{l([\alpha])(-Q(\omega_{i})-(\rho(K)-\epsilon))}\\
&\leq& C e^{t\left ( 
-Q(\omega_{i})-(\rho(K)-\epsilon)\right )}
\sum_{[\alpha] \in [\Gamma_{hyp}](t):
\delta_{[\alpha]} \in K_{i}}  e^{\int_{[\alpha]} (f+\omega_{i})}.
\end{eqnarray*}
For $t$ sufficiently large, it follows from Theorem $1$ $(1)$,
\begin{eqnarray*}
\nu_{t}(K)& \leq& 
C\sum_{i=1}^{l}e^{t(P(f+\omega_{i})+\epsilon)}e^{-t(P(f)-\epsilon)}
e^{t(-Q(\omega_{i})-(\rho(K)-\epsilon)}\\
&=&Cl e^{t(-\rho(K)+3\epsilon)}.
\end{eqnarray*}
Take the logarithme, divide by $t$ and take the $\limsup$,
\[
\limsup_{t\rightarrow \infty}\frac{1}{t}\log \nu_{t}(K) \leq
-\rho(K)+3\epsilon. 
\]
$\epsilon$ being arbitrary, this proves Lemma $3$.

\end{document}